Exploring Gender Differences in Tertiary Mathematics-Intensive Fields: A Critical Review of Social Cognitive Career Theory

| Huayu Gao | Tanya Evans | Gavin Brown |
| University of Auckland | University of Auckland | University of Auckland |

*Social Cognitive Career Theory (SCCT) has been extensively employed to elucidate the enduring gender differences in mathematics-intensive fields, with a particular emphasis on the complex interplay of motivational factors and extra-personal influences contributing to the underrepresentation of women. Although a plethora of empirical studies corroborate SCCT, three crucial aspects for refinement have come to the fore. First, the theory should place a more substantial emphasis on how cultural and contextual diversity influences academic choices. Second, given the dynamic nature of motivation, which evolves over time, more longitudinal analyses are imperative to capture their temporal trajectory, in contrast to the predominantly cross-sectional empirical studies. Finally, considering the intricate interplay between emotion and motivation, integrating the dimension of emotion into SCCT would significantly augment its explanatory power and provide a more comprehensive understanding of academic selection processes.*

*Keywords:* Gender Differences; Mathematics-intensive Academic Fields; Social Cognitive Career Theory; Women in STEM

  Gender disparities in tertiary mathematics-intensive fields and subsequent career trajectories remain a pervasive issue in most Western countries, with women consistently underrepresented. This persistent challenge is particularly evident in tertiary education, where women's participation in disciplines such as geoscience, engineering, economics, mathematics, computer science, and the physical sciences (chemistry and physics) is notably lower compared to their male counterparts (Ceci et al., 2014). Alarmingly, women's representation in computer science and electrical and electronic engineering has not only stagnated but has declined over time in Spain (Olmedo-Torre et al., 2018). Moreover, the gender gap transcends the boundaries of academia and permeates the workforce, with women's participation in math-intensive careers ranging from a mere 14% to 24% across various countries (Piloto, 2023), further underscoring the systemic nature of this disparity.

  The persistent gender disparity in mathematics-intensive domains within tertiary education and its far-reaching implications for women's career trajectories necessitate a comprehensive exploration of the underlying factors perpetuating this inequality. Specifically, the underrepresentation of women in math-intensive fields serves as a critical antecedent to their subsequent underrepresentation in related career paths, exerting a profound influence on various aspects of society, including economic growth, wage equality, and the effective utilization of women's potential (Card & Payne, 2021). This underutilization of women's potential is particularly alarming, given that girls often outperform boys at most educational levels (OECD, 2015). To effectively address this multifaceted issue, Social Cognitive Career Theory (SCCT) has emerged as a valuable framework for elucidating the persistent gender gap, providing nuanced insights into the complex interplay of individual, social, and institutional factors that influence academic selection and career trajectories.

  SCCT offers a robust framework for understanding academic choices and career development among boys and girls; however, empirical evidence indicates that further refinement is necessary to enhance the theory's explanatory and predictive capabilities. This paper critically elucidates the limitations of SCCT, particularly its applicability across diverse cultural contexts, its reliance on cross-sectional empirical research despite the theory's dynamic nature, and its underestimation of emotional factors in academic decision-making. By identifying these limitations and proposing potential improvements, this theoretical paper aims to contribute to the ongoing evolution of SCCT, ensuring that it remains a relevant and effective tool for understanding individuals' choices and performance among boys and girls in increasingly complex and diverse environments.

## Social Cognitive Career Theory: Key Concepts and Relations

Social Cognitive Career Theory, developed by Lent and Brown (1994), is a comprehensive framework that elucidates the intricate processes involved in academic selection and career development. Rooted in Bandura's (1986) social cognitive theory, SCCT emphasizes the interplay between individual agency and environmental factors in shaping academic choices and career trajectories. The theory posits three interconnected models—interest, choice, and performance, each centered around the core constructs of self-efficacy, outcome expectations, and goal representations. These constructs interact dynamically with personal factors, contextual influences, and learning experiences to shape career-related behaviors and decisions.

## The Distinctive Roles of Self-Efficacy and Outcome Expectations in SCCT

Self-efficacy and outcome expectations are two key constructs in Bandura's (1986) social cognitive theory that play crucial roles in shaping individuals' selection processes. Self-efficacy refers to individuals' appraisals of their capability to organize and execute the requisite actions to attain desired performance outcomes. Crucially, self-efficacy is distinct from objectively assessed abilities. Social cognitive theory posits that human capabilities are malleable rather than fixed, asserting that successful performance in complex or demanding tasks typically necessitates not only foundational skills but also a robust sense of efficacy to effectively deploy those skills (Bandura, 1991). While self-efficacy focuses on individuals' confidence in their ability to execute actions, outcome expectations pertain to the anticipated consequences of those actions. Bandura (1986) classified outcome expectations into three categories: physical (e.g., tangible rewards), social (e.g., external validation), and self-evaluative (e.g., self-satisfaction). These expectations can exert a significant influence on career-related decisions and behaviors.

Bandura (1986) acknowledged the dual role of self-efficacy and outcome expectations, but he contended that these two forms of belief often differ in their relative potency, with self-efficacy serving as a more influential determinant of behavior. To illustrate this point, there are numerous instances in which individuals may anticipate valued outcomes resulting from a particular course of action but refrain from pursuing such action if they harbor doubts about their capabilities. However, in situations where outcomes are only tenuously linked to performance quality, outcome expectations may contribute to motivation and behavior independently (Bandura, 1989). This scenario may be especially pertinent to career development, as the intricacies of academic and career environments often yield only imperfect connections between performance quality and outcomes. For instance, an individual with a high ability for mathematics might opt to eschew mathematics-intensive career fields if they lack the support of significant family members; their anticipated outcome is far away from math.

**The Structure of SCCT: Interest, Choice, and Performance Model**

*Interest model.* Self-efficacy and outcome expectations, which incorporate values, are crucial in forming interests. Through repeated activity engagement, modeling, and feedback, children and adolescents refine skills, develop performance standards, and establish self-efficacy in specific tasks while acquiring outcome expectations (Lent & Brown, 1994). The development of self-efficacy and outcome expectations contributes to the formation of interests both directly and indirectly. These emerging interests lead to intentions for further activity exposure, increasing the likelihood of task selection and practice. Activity involvement results in performance attainments, such as successes and failures, which prompt the revision of self-efficacy and outcome expectancy estimates (Lent & Brown, 1994).

*Choice model.* The choice model builds upon the foundation of the basic interest model, extending it to encompass the developmental aspects of academic or career decision-making. This model conceptually divides the choice phase into several interconnected processes (Lent & Brown, 1994). Firstly, individuals express a primary choice goal, selected from their major interests. This choice goal represents the intention to pursue a specific career path or engage in a series of actions, such as aspiring to become an engineer. Secondly, individuals take concrete steps to implement their choice, which may include enrolling in relevant training programs or selecting an appropriate academic major. Finally, the model incorporates a feedback loop, wherein subsequent performance attainments, such as academic successes or failures and admission outcomes, shape future career behavior and choices.

*Performance model.* The performance model presented here adopts a broad definition of performance, encompassing both the level of achievement (e.g., academic grades) and indices of behavioral persistence (e.g., stability in a chosen major). This model emphasizes the significant roles played by ability, self-efficacy, outcome expectations, and performance goals while omitting interest as an intermediate mechanism. Lent & Brown (1994) posit that interest is more closely linked to selecting career or academic activities rather than setting performance goals. As the choice model suggests, choice goals contribute to determining the performance domains (e.g., work tasks) that an individual will pursue; however, the quality of performance attained may be partially contingent upon the level of one's performance goals. For instance, students who opt for an engineering major (choice goal) also establish goals pertaining to their grade performance in the requisite courses. Aspiring to achieve an A in math (performance goal) serves to regulate one's subsequent course-related behaviors (e.g., study time allocation).

**Sociocultural Influences on Gender Disparities in Career Choice and Performance**

The SCCT provides a robust framework for understanding how sociocultural factors contribute to gender disparities in career choice and performance from a motivational perspective. Specifically, the theory emphasizes the influence of personal factors, such as ethnicity and gender; contextual factors, including cultural norms, gender role expectations, and emotional support; and experiential factors, like personal achievements, vicarious learning, social persuasion, and physiological states, on career-related outcomes (Lent & Brown, 1994). These factors serve as antecedents of motivational variables, moderators of key relationships, or direct facilitators or barriers to career development. Furthermore, SCCT posits that gender differences in motivation are primarily influenced by social contexts. While acknowledging the role of innate aptitudes and interests, the theory suggests that these predispositions can be molded by environmental factors throughout an individual's life, highlighting the significance of the social contextual factors in shaping career-related outcomes and contributing to gender disparities.

# Social Cognitive Career Theory: Empirical Evidence on Gender Perspectives

**Gender Differences in Motivation**

A growing body of research consistently demonstrates the predictive validity of self-efficacy and outcome expectations in relation to interests, choices, and performance, particularly in mathematics-intensive domains (Lent et al., 2017; Sheu et al., 2010). These findings underscore the central role of motivational variables in shaping individuals' math-related behaviors and outcomes. However, research also indicates that these motivational variables do not operate in a vacuum; they are influenced by broader socio-cultural factors, including gender-role identities. For instance, due to societal expectations and the internalization of gender norms, girls often exhibit lower self-efficacy and outcome expectations, as well as reduced interest in math-intensive fields compared to their male counterparts (Luo et al., 2021; Lv et al., 2022).

**The Impact of Family and School on Gender-Based Motivational Differences**

As essential carriers of socio-cultural influences, family and school have been widely recognized as critical determinants of choices and performance. The family, as the primary social unit, provides children with their initial exposure to gender-specific expectations and roles. This early socialization within the family context lays the foundation for the development of gender-specific beliefs and behaviors. Schools, on the other hand, provide opportunities to refine these beliefs and behaviors. Therefore, the complementary relation between family and school is crucial in the evolution of children's career aspirations.

*Family factors*. The family factors, including socioeconomic status (SES), parental beliefs, and parental behaviors, exert a profound influence on individuals' academic achievements and career aspirations. SES, a multidimensional construct encompassing income, education, and occupation, has been consistently linked to children's academic choices. Gao et al. (2024) demonstrated that SES accounted for a substantial 50% of the variance in adolescents' mathematics achievement in New Zealand. Mathematics achievement serves as a critical prerequisite for students' selection of math-related majors (Wang & Degol, 2013). Higher SES families tend to provide greater resources, support, and extracurricular activities that foster the development of mathematical skills (Bandura & Bussey, 2004; Wang & Degol, 2013), thereby fostering the formation of interests and influencing future choices. Conversely, low-SES families often prioritize meeting their children's fundamental daily necessities over fostering the pursuit of mathematics-related competencies, with girls' academic performance and development receiving even less attention and support.

Beyond SES, parental beliefs about gender also exert an important influence on children's confidence and career aspirations. Empirical evidence suggests that parents' perceptions of their children's abilities can have a more significant impact on their offspring's self-beliefs and career goals than actual performance (Fredricks & Eccles, 2002; Parsons et al., 1982). More specifically, parents who endorse gender stereotypes in math tend to underestimate their daughters' abilities while overestimating their sons' abilities, and they may even place less importance on their daughters' participation in tertiary math-intensive majors (Parsons et al., 1982; Nicoletti et al., 2022). Moreover, parental behaviors often vary based on gender, inadvertently reinforcing gender disparities. Research has revealed that parents are more inclined to engage in scientific discussions with their sons and provide them with related resources while unintentionally limiting opportunities for their daughters (Crowley et al., 2001; Meece et al., 2009). This differential treatment can perpetuate gender stereotypes and contribute to the underrepresentation of women in math-intensive areas (Jacobs et al., 2005).

*School Factors.* Alongside family influences, schools and tertiary institutions also play a crucial role in shaping gender-specific beliefs and influencing students' motivation and career aspirations. Professors, teachers and peers are key agents in this process, with their beliefs, expectations, and interactions contributing to gender disparities in math-related participation. Specifically, teachers' beliefs and expectations about gender-specific stereotypes can significantly impact students' motivation. For example, teachers may emphasize girls' efforts over abilities in mathematics, leading to lower expectations and diminished confidence among girls (Patrick et al., 2007; Meece et al., 2006). However, the effect of teacher support can vary by grade level, with elementary school environments often favoring girls, while secondary schools tend to emphasize performance goals that may disadvantage them (Wang & Degol, 2013).

Peers also play a significant role in influencing math-related achievement and participation. Girls in supportive peer groups tend to excel in mathematics and science; however, the Big-Fish-Little-Pond Effect suggests that while being surrounded by high-achieving peers can provide academic stimulation and motivation, it may also lead to a diminished self-concept. This phenomenon occurs because individuals tend to compare their performance to that of their peers, and in a group of high achievers, even objectively strong performers may perceive themselves as less competent, thus lowering their academic self-concept (Leaper et al., 2012; Marsh, 1984). In addition to the influence of teachers and peers, it is important to consider the potential impact of single-sex schooling on girls. Docherty et al. (2020) discovered a peculiar discrepancy in engineering enrolments at the University of Canterbury (New Zealand) between 2005 and 2017. The study revealed that the proportion of female enrolments in engineering from single-sex high schools was disproportionately high, accounting for 56% of all female students in engineering, despite the fact that only a small proportion of girls (15.73%) attend single-sex schools (MoE, 2023). This finding suggests that single-sex settings may play a role in shaping girls' academic choices. Moreover, Smith and Evans (2024) reported that girls attending single-sex schools perform better in mathematics and science than their co-educational counterparts in New Zealand. They suggested that single-sex environments may reduce stereotype threat for girls by providing them with a more supportive and nurturing learning context, thereby fostering their confidence and interest in traditionally male-dominated fields. However, the research on the effectiveness of single-sex schooling is complicated, with some studies suggesting that single-sex schools do not have any significant advantages for girls in pursuing math-related goals compared to co-educational settings (Gao et al., 2024).

**Social Cognitive Career Theory: Limitations and Pathways for Future Development**
The aforementioned studies offer considerable empirical evidence that substantiates the SCCT, especially concerning its capacity to elucidate gender disparities in mathematics-intensive domains. Nevertheless, while these studies corroborate the SCCT framework, it is imperative to acknowledge that certain research also underscores the theory's limitations. These constraints encompass the inadequate emphasis on cultural diversity, the predominant reliance on cross-sectional research, and the absence of an emotional dimension. The limitations identified through our summary provide valuable guidance for future scholarly efforts to refine and expand the SCCT framework, thereby strengthening its explanatory capacity and applicability across diverse contexts.

**Placing Stronger Emphasis on Cultural Diversity to Strengthen SCCT**
The SCCT underscores the crucial role of sociocultural factors in shaping gender differences, but it falls short in providing a comprehensive analysis of the social environment due to its Western origins, which may limit its universal applicability across diverse global contexts. Despite distinguishing between the social and biological concepts of gender and ethnicity, with a primary emphasis on the former, the SCCT has not adequately addressed the influence of cultural backgrounds in different countries on individuals' motivation and career choices. Recent cross-cultural studies investigating the application of SCCT have yielded inconsistent results, underscoring the complex interplay between cultural factors and academic choices. For example, Wang et al. (2023) identified significant gender disparities in the application of SCCT within the Chinese educational system, with girls exhibiting lower interest and requiring greater confidence in mathematics compared to boys. Conversely, Gao et al. (2024) found no such gender discrepancies in New Zealand. Moreover, Michaelides et al. (2019) observed that in some Eastern European and Islamic countries, girls outnumbered boys in math-intensive fields and demonstrated higher confidence in mathematics. These divergent findings highlight the need for a more nuanced understanding of the cultural factors influencing academic selection processes and suggest that the explanatory power of SCCT may vary considerably across different cultural contexts.

To address the limitations of SCCT in capturing the complexities of cultural influences on academic and career decision-making processes, future research should prioritize the investigation of cultural background factors and their underlying mechanisms in shaping individuals' educational and vocational trajectories. This endeavor necessitates the development of a substantial body of cross-cultural empirical research alongside the creation of culturally specific models that can effectively capture the nuances of diverse contexts. For instance, in Iran, where women's status is relatively low, science, technology, engineering, and mathematics (STEM) fields are often perceived as pathways to higher income and social status, potentially leading women to choose these fields to improve their social standing (Stoet & Geary, 2018; Michaelides et al., 2019). In contrast, in some Western countries with higher levels of gender equality, girls typically excel across various disciplines and face a wider range of career options. This diversity of choices may result in a lower proportion of girls choosing STEM fields, as they are more likely to make decisions based on their absolute abilities, whereas boys are more likely to make choices based on their relative abilities (Wang & Degol, 2013). As evidence gradually accumulates, it will serve as a foundation for refining and expanding the theory, enabling it to more accurately account for the intricacies of academic selection processes across various cultural landscapes.

**Employing Longitudinal Studies to Validate SCCT**

The cross-sectional research designs employed by the vast majority of empirical investigations examining the tenets of SCCT pose inherent limitations in their ability to fully capture the dynamic feedback mechanisms central to the theory. Specifically, SCCT posits a complex interplay of feedback loops in which an individual's self-efficacy beliefs and outcome expectations exert a profound influence on their choices and performance. In turn, the consequences of these choices and the resulting achievements serve to reshape future self-efficacy beliefs and outcome expectations, creating a continuous cycle of reciprocal determinism (Lent et al., 1994). However, cross-sectional studies are unable to adequately test these intricate relations of reciprocal determinism.

To fully elucidate the intricate dynamics posited by SCCT, a longitudinal approach is indispensable. Such research designs facilitate the exploration of temporal relations and enable the mapping of individual trajectories across extended timeframes (Ployhart & Vandenberg, 2010). By gathering data at multiple junctures, researchers can delve into the nuanced ways in which self-efficacy and outcome expectations evolve in response to diverse experiences and interventions, thereby offering a more comprehensive and ecologically valid assessment of SCCT (Navarro et al., 2014). Furthermore, longitudinal studies possess the unique capacity to illuminate the directionality of relations between variables, tackling questions of causality that remain elusive in cross-sectional designs (Zapf et al., 1996).

**Integrating Emotional Dimensions to Broaden SCCT**

The SCCT has a notable limitation in its insufficient emphasis on the role of emotions in shaping individuals' academic choices. While the theory primarily focuses on explaining these choices from a motivational perspective, highlighting the interplay between self-efficacy, outcome expectations, and interest, it does not adequately address the significant influence of emotions on cognitive functions crucial to academic decision-making. If the value of outcome expectation negatively influences individuals' math achievement, it is challenging to explain using SCCT. However, by incorporating emotion into the SCCT framework, the result becomes plausible. It can be argued that a higher value may induce math anxiety, consequently lowering achievement. More specifically, emotions play a critical role in shaping cognitive processes such as attention, memory, and problem-solving (Lerner et al., 2015), which are essential components of academic-related decisions. Positive emotions have been associated with higher levels of engagement, deeper learning, and improved performance, contributing to stronger self-efficacy and more favorable outcome expectations (Reschly et al., 2008). In contrast, negative emotions can disrupt these cognitive processes, leading to decreased motivation and poorer academic outcomes (Pekrun, 2006). Consequently, the lack of sufficient attention to the emotional dimension within the SCCT framework may limit its explanatory power in fully capturing the complex dynamics of the academic choice process.

To address this limitation, future research should focus on integrating emotional dimensions into the SCCT framework. For example, Pekrun's (2006) Control-Value Theory (CVT) combines motivation and emotion to predict individuals' choices and performance in educational contexts. In CVT, control refers to individuals' beliefs about their abilities, i.e., whether they believe they can complete a task, while value includes intrinsic value (such as interest) and extrinsic value (such as the importance of the task). When individuals perceive their abilities as insufficient (low control) but simultaneously consider the task outcomes to be crucial for them (high extrinsic value), they often experience negative emotions such as anxiety or helplessness, which can adversely affect their academic performance. Therefore, by incorporating emotional factors into the SCCT framework, researchers can more comprehensively explain the dynamic interactions between emotions, cognition, and situations, thus establishing a more realistic theoretical model to describe individuals' academic choices and behaviors. This integration not only expands the theoretical perspective of SCCT but also significantly enhances its applicability to diverse populations and contexts, especially those groups and environments where emotional responses vary significantly.